# Efficient, Fast, and Fair Voting Through Dynamic Resource Allocation in a Secure Election Physical Intranet

Tiankuo Zhang, Benoit Montreuil, Ali V Barenji, Praveen Muthukrishnan
Physical Internet Center, Supply Chain & Logistics Institute
School of Industrial and Systems Engineering
Georgia Institute of Technology, Atlanta, USA
Corresponding author: tzhang417@gatech.edu

*Abstract:* Resource allocations in an election system, often with hundreds of polling locations over a territory such as a county, with the aim that voters receive fair and efficient services, is a challenging problem, as election resources are limited and the number of expected voters can be highly volatile through the voting period. This paper develops two propositions to ensure efficiency, fairness, resilience, and security. The first is to leverage Physical Internet (PI) principles, notably setting up a "secure election physical intranet" (SEPI) based on open resource sharing and flow consolidation between election facilities in the territory. The second is to adopt a smart dynamic resource allocation methodology within the SEPI based on queueing networks and lexicographic optimization. A queueing model is developed to provide feasible combinations of resources and individual performances for each polling location by considering layout and utilization constraints. A two-stage lexicographic optimizer receives the queueing model's outputs and finds an optimal solution that is less expensive, fast, and fair. A scenario-based case study validates the proposed methodology based on data from the 2020 US Presidential Election in Fulton County, Georgia, USA.

*Physical Internet (PI) Roadmap Fitness*: ☒ PI Networks, ☒ Access and Adoption

*Targeted Delivery Mode-s*: ☒ Paper, ☐ Poster, ☐ Flash Video, ☒ In-Person presentation

## 1  Introduction

Election is paramount to any democratic framework worldwide, yet long, uneven, and volatile waiting lines threaten this critical cornerstone. Within the United States, voters and election system administrators have noticed a rising waiting time in an increasing number of polling locations. During the U.S. in 2012, many polling locations in battleground states saw voters waiting for hours on election day (Famighetti et al., 2014). These long waiting lines can discourage people from voting and induce a huge national economic cost, as even a modest wait of 10-15 minutes from each voter can cost the U.S. approximately $500 million (Stewart III et al., 2013). Additionally, long waiting times disproportionately affect different polling locations, which can lead to unfairly treated voters' frustrations and potential impacts on election results. For example, during the U.S. Presidential Election of 2020, certain polling locations in Florida reported average waiting times of over 80 minutes, while other polling locations nearby had less than 10-minute waiting times (Barenji et al., 2023). Other studies also suggest that certain voter populations, such as racial minorities and lower-income communities, are more likely to experience long waiting times (Stewart III, 2013; Pettigrew, 2017). Addressing these concerns is imperative to safeguarding the democratic process and ensuring

fair electoral participation for all citizens, and one pivotal solution is effective election resource allocations with fast and fair services.

Key resources to be allocated in an election system include poll pads, ballot marking devices (BMDs), and scanners. The allocation is often challenging due to voter turnout volatilities. Many factors can affect voter turnouts, drastically differing between county precincts and even between the same precinct's hours (Matsusaka et al., 1999). Resource limitations, location constraints, and governmental regulations also bound allocation plans. For instance, based on the Centers for Disease Control and Prevention's (CDC) recommendations, the social distancing regulation for COVID-19 greatly affected resource allocation during the U.S. Presidential Election of 2020 and led to many polling locations' shortness in providing efficient services (Sullivan, 2020). Past decision-making in polling resources mostly involved fixed population-based plans, in which all polling locations would receive unchanging amounts of various resources at the beginning of early elections, and the amounts usually relied on registered voter populations, historical plans, or the busiest anticipated day (Edelstein, 2006). However, fixed apportionment can lead to a waste of resources and cannot respond to a surprising uprising of voters. For instance, throughout the early-election period of the U.S. presidential election in 2020, many polling locations in Georgia experienced hours of waiting due to the number of early-election voters being higher than expected based on the 2016 turnout data.

In this paper, we develop two propositions to ensure efficiency, fairness, resilience, and security. First, we leverage Physical Internet (PI) principles (Montreuil, 2011), notably setting up a Secure Election Physical Intranet (SEPI) based on open resource sharing and flow consolidation between election facilities in the territory. The SEPI incorporates three considerations and ensures a secure and resilient allocation process. Second, we adopt a smart dynamic resource allocation methodology based on queueing network analysis and lexicographic optimization.

We assess the value of our propositions through comparisons with the fixed resource allocation plan derived from Barenji et al. (2023) for the 2020 US Presidential Election in Fulton County, Georgia, USA. Our approach outperforms the fixed official resource allocation plan regarding total resource and cost requirements, resource utilization, and voters' estimated robust wait time and fairness. We also identify critical polling locations that require policymakers' attention.

The full paper is organized as follows. Section 2 presents the related literature. Section 3 demonstrates the SEPI. Section 4 proposes the dynamic resource allocation framework. Section 5 analyzes the methodology's performance through a case study of three scenarios. Section 6 summarizes the contributions, limitations, and future work directions.

## 2 Related Literature

Researchers have conducted many related studies on election resource allocations and their assessments. Allen and Bernshteyn (2006) suggested using queueing network models in estimating and mitigating voter waiting times on election day. Its case study proved its effectiveness in reducing voters' average voting time. Allen et al. (2020) explored the application of indifference zones in determining whether a set of allocated resources can guarantee acceptable wait times in a polling location. A general binary search algorithm is also employed to improve the run time when seeking the optimal allocation combination. However, this set of procedures only yields limited options with high performances, and satisfying all these options for all polling locations from a managerial perspective can be challenging as election resources are often limited. Furthermore, Allend and Benshteyn's (2006) queueing network model and Allen et al.'s (2020) indifference-zone binary search approach only explored fixed resource allocation plans. They could not respond to demand fluctuations, especially when a polling location faces an unexpected demand rise.

Previous PI container studies have introduced and proven their specialty in improving security, efficiency, and reliabilities (Montreuil, 2011; Montreuil et al., 2016; Sallez et al., 2016). Their most fundamental pillar in protecting contained objects made them suitable for applying election systems, as election resources are extremely sensitive to the outside world's contamination. PI containers' modularity and traceability can also guarantee efficiency and reliability during resource transfers, which is critical in election systems, ensuring all resources are delivered quickly and correctly. Furthermore, PI containers' specialties in sustainability can help alleviate economic and environmental stresses running the election system.

Another related topic in election resource allocations is voter turnout predictions. Two major categories of research are conducted in this field: individual-level analysis and population-level analysis. Individual-level analyses predict whether an individual voter will vote based on his/her utility cost function, of which positivity indicates whether a voter will vote (Lacy et al., 1999; Garcia-Rodriguez, 2020). Population-level analyses predict voter turnout rates of the entire population of a county, state, or nation (Huang, 2016; Turner, 2024). Specifically, Huang (2016) utilized Dynamic Learn Models (DLMs) to predict voter turnouts based on historical data.

## 3   Proposed Secure Election Physical Intranet

Security, ensuring all machines are functional and not tampered with, is critical to election resource allocations. While a piece of machine is transferred, it is exposed to more damage or interference from maleficent individuals or groups. Therefore, we aim to develop procedures to enhance machine security during and after transfers.

To protect resources from being damaged or touched during transfers, our methodology leverages Physical Internet concepts and principles to be implemented within a secured physical intranet for the election territory, with open resource sharing and flow consolidation among polling locations, election hubs, and warehouses (Figure 1). Resources are to be stored and moved in sealable PI containers with modular sizes and connectors such that transferred machines in compact and fitting sizes, greatly enhancing the security of the entire system (Montreuil et al., 2016). We also suggest only using localized machine transfers such that only polling locations of a certain managerial range, such as the commission district illustrated in Figure 1, can transfer resources to each other. By limiting the resource transfer range, machines' travel distances are greatly reduced, reducing exposure to dangers and the chance of being interfered with. Moreover, polling locations within a managerial range will not have to cross administrative barriers to communicate with other election management, improving dynamic resource allocation's speed and safety. After election resources are safely transferred to a new polling location, we encourage staff to examine the transferred machine fully. The examination protocols should be extensive and standardized so that no damaged machines would be put into actual use on the next day.

The above considerations are realized within our cost calculator model and the cost functions within the optimization models. The cost calculator model computes the cost matrix, in which polling locations from two different managerial ranges have a transfer cost of infinity that directly prevents optimization models from selecting such options. The optimization model's cost function captures the costs of differently sized PI modules, transfers, and examinations.

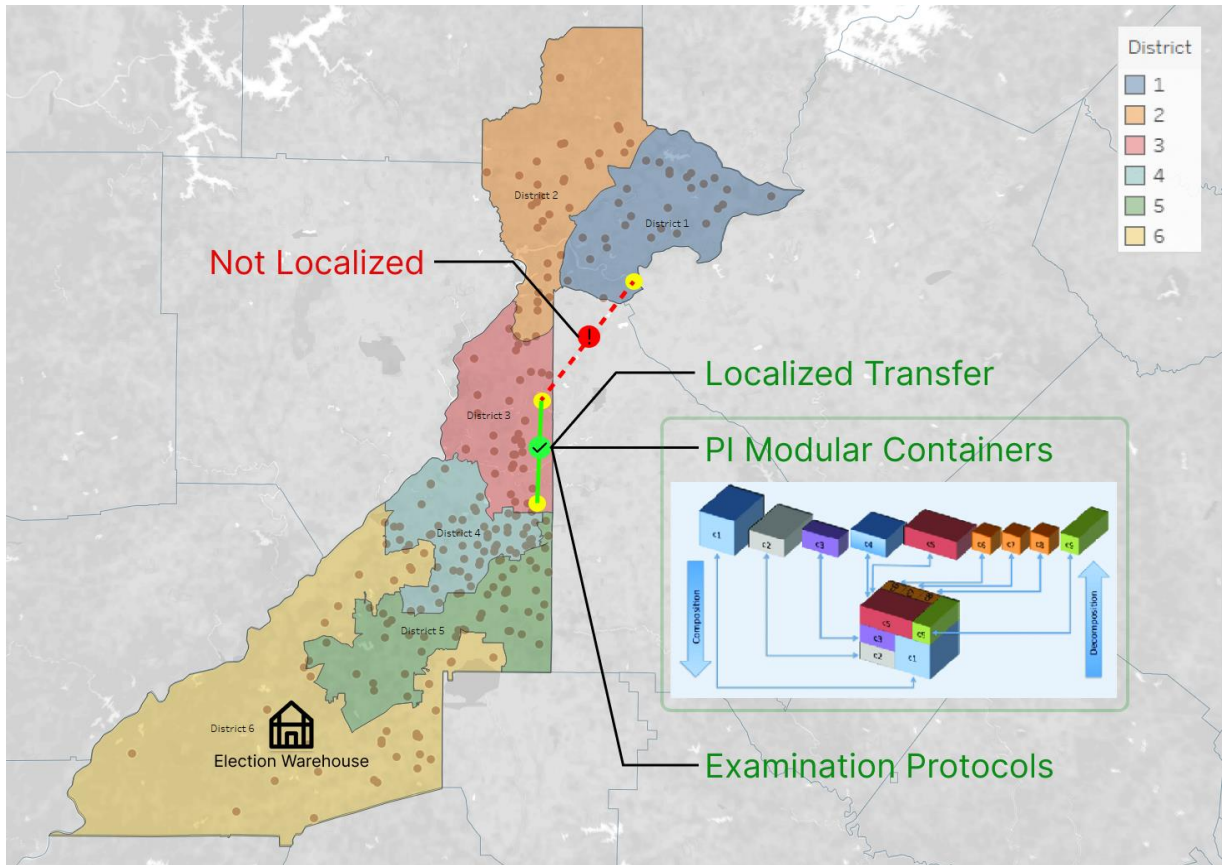

*Figure 1: Illustrating the Proposed SEPI in Fulton County, GA (Montreuil et al., 2011)*

In this context, the paper proposes a dynamic resource allocation methodology for supporting resource allocation and transfer of resources to polling locations across the SEPI territory, such as a county. For each polling location throughout the early voting period (e.g., 19 days in the U.S. in 2020), it monitors voter hourly turnout, predicts voter turnout in the upcoming election days, and seeks an optimal resource allocation and transfer plan involving all polling locations. Based on the transfer plan, each polling location, including the territory's election resource warehouse, transfers a designated amount of resources of each type to each other during nights when polling locations are closed to the public, abiding by election regulations.

## 4   Proposed Dynamic Resource Allocation Framework

Figure 2 illustrates the proposed dynamic resource allocation methodology's framework, which contains one major data pool and five computation models: the cost calculator, demand model, queueing network model, and a lexicographic optimization system composed of a first-stage optimizer and a second-stage optimizer. The cost calculator takes polling location latitude and longitude information as inputs and computes a cost matrix, of which each entry contains the transfer cost between two polling locations. The demand model intakes historical turnout and turnout data from previous early election days of the current election. It utilizes linear regression and statistical forecasting to predict hourly arrival for each polling location each day until the election ends. To improve the overall system's performance and consider the up-to-date demand fluctuation, our proposed demand model, queueing network model, and optimization models should be run at the end of every pre-election day for the most up-to-date plan. Throughout the early election period, the election officials are to direct the realization of optimized inter-location transfers with a designated amount of resources of each type during non-voting hours based on the latest refreshed plan.

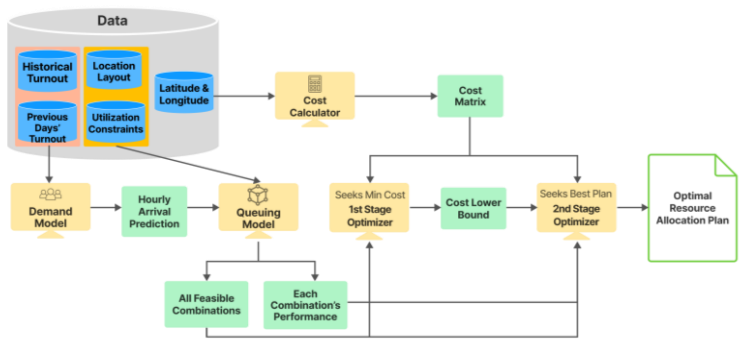

*Figure 2: Illustrating the Proposed Dynamic Resource Allocation Framework*

## 4.1 Queueing Network Analysis Model

The queueing network model takes the demand model's hourly arrival predictions, location layout constraints, and utilization rate constraints as inputs and computes all feasible resource combinations and each combination's performance for each polling location for all upcoming days until the election ends.

For each polling location on each day, we first select all feasible resource combinations that satisfy the location's layout constraint. Based on each feasible resource combination, the currently modeled queueing network considers three queues that simulate a full voting process: check-in at poll pads stations, casting votes at BMDs stations, and the scanner station, as displayed in Figure 3. A performance array is also initialized to record the performance values. The queueing network is then run for many iterations, such that results are stable and rigorous according to the Central Limit Theorem. Within each iteration, a voter-arrival array is built based on the hourly arrival prediction and used to compute and store the resource combination's performances, such as waiting time and each machine's utilization rate. Upon completing all iterations, we evaluate the current combination's utilization rates based on the stored performance values. If the utilization rates do not meet the utilization constraints, we will abandon the current combination. Otherwise, we continue to compute its robust waiting time and store it for future usage.

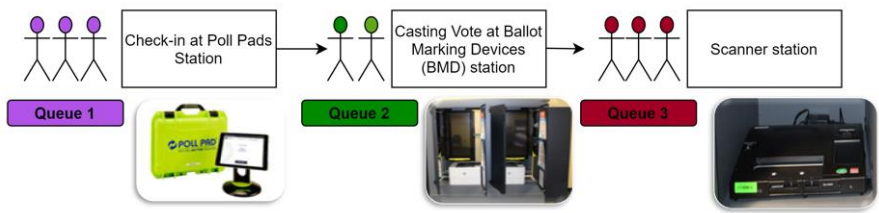

*Figure 3: Demonstrating Queues within a Queueing Network*

## 4.2 Allocation Optimization Model

The optimization model selects resource combinations from the queueing network model's computation results. The selected resource combinations are satisfied by utilizing a set of combination constraints, and the optimization model determines which location should transfer how many of which resources to other specific locations. A modified indifference-zone approach categorizes each resource combination's performance by translating its achieved robust waiting time (e.g., at 99.7%) into a performance score from 0 to 1, which does not need to exhibit a linear relationship and can be tuned as pertinent. This application results in a more effective resource allocation plan by preventing the optimization model from allocating more resources to polling locations that already satisfy most voters' expectations. For example, the Presidential Commission on Election Administration (PCEA) suggests that most voters are willing to accept waiting times less than 30 minutes. As election resources are limited, in this illustrative situation, if a certain polling location's resource combination enables a 99.7% confidence waiting time of fewer than 30 minutes, an efficient resource allocation plan may not allocate more resources to this location by creating an indifference zone of all wait times below 30 minutes.

Our dynamic election resource allocation methodology aims to find less expensive, fast, and fair allocation plans, translating into three optimization objectives: minimize the total cost, maximize the total polling location performance score, and minimize the performance score gap. The third optimization goal would require searches for the best and the worst performances, resulting in a nonlinear optimization model that requires a long runtime. To resolve this challenge, we employ a pair of lower and higher performance score bounds within the optimization model, such that all selected resource combinations for all locations must have performance scores larger than the lower bound and smaller than the higher bound. By tightening the performance score bounds, we can create a set of fairness constraints to ensure a small performance gap between all polling locations on all days without nonlinear optimization models.

Cost minimization and performance score maximization conflict, as lowering the total allocated resources and limiting the number of transfers unavoidably result in certain polling locations not receiving enough resources for the best performance score. Therefore, our methodology employs lexicographic optimization in solving for the multi-objective system. Specifically, we rank cost minimization as the first objective and performance score maximization as the second and proceed according to the two-stage optimization methodology in Figure 4.

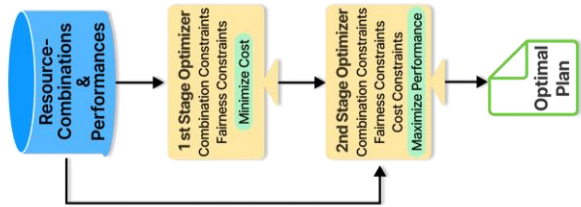

*Figure 4: Showing the Optimization Model's Procedures*

The stage-one optimizer takes the current resource allocation plan, combinations generated by the queuing model, and the cost matrix as inputs and finds the lowest possible total cost among all possible plans, satisfying the robust waiting time and fairness constraints. The stage-two optimizer utilizes all of the stage-one optimizer's inputs and constraints, and the relaxed cost lower bound produced by the stage-one optimizer as an additional constraint and finds the optimal plan that ensures low cost, low wait, and high fairness across all polling locations throughout the election.

## 5 Case Study

Fulton County, Georgia, is the state's most populous county, with a population of over 1 million, according to the U.S. Census Bureau in 2022. During the 2020 presidential election, Fulton County defined 238 polling locations, with 149 supporting early voting of 19 days. The turnout rate was approximately 65%, and around 60% of voters chose to vote during early elections. Based on Fulton County's 2020 presidential election, we construct two demand-fluctuation scenarios: a high early voting scenario and a low early voting scenario.

In both scenarios, we assume all polling locations each day exhibit the same hourly arrival percentage of their estimated voter turnouts as their counterparts on election day. The ratio between the daily voters of each polling location and the total voters of all polling locations on each day remains fixed. The total number of registered voters, the total turnout, and the number of absentee voters (voters by mail) are all consistent with Fulton County's 2020 presidential election. We realize our previous applications of PI modules and examinations within the cost functions of optimization models. Furthermore, we assume each polling location only transfers resources to polling locations within the same commission district, corresponding to our localized transfer concern, and there are six commission districts.

The high early voting scenario assumes that 75% of the voters choose early voting. As the number of absentee voters remains unchanged, fewer voters vote on election day. Similarly, the low early voting scenario assumes 45% choose early voting, which results in more voters voting on election day. We made assumptions about these two scenarios based on the observed trend from the U.S. presidential election from 2016 to 2020. In 2016, approximately 45% of voters chose early voting, increasing to 60% in 2020. However, considering COVID-19's impact and the CDC's recommendations, it is unclear whether this trend will continue (Santana et al., 2020). Thus, we create two scenarios: the high early voting scenario corresponds to a continuation, and the low early voting scenario reflects a return.

We compare our allocation plan's performance to the fixed allocation plan derived from Barenji et al., 2023. Although fixed, this allocation plan is developed through a multi-agent-based simulation platform and is superior to the official allocation (Barenji et al., 2023). We use a queueing network model separate from the one applied in the methodology to examine each polling location's performances during the early voting period and election day. We split wait times into 6 ranges and utilize pie charts to demonstrate the percentage of polling locations within each range during early election and election day. Through the comparisons, polling locations that require more resources or expansion are identified.

Figure 5's left side compares the fixed and dynamic allocation plans' 99.7% rigorous wait times in the high early voting scenario. The fixed allocation plan allows all polling locations to wait less than 30 minutes throughout the early election. However, around 8% of polling locations saw voters waiting more than 30 minutes on the election day, within which most polling locations' wait time lies in the 30-45 minutes range, and multiple polling locations' voters are estimated to wait more than 150 minutes. Compared to the fixed plan, dynamic resource allocation enables all voters at all polling locations to vote after waiting for less than 30 minutes during early voting and election day, indicating a fairer waiting time distribution.

Figure 5's right side compares the fixed and dynamic allocation plans in the low early voting scenario. The fixed allocation plan continues to provide efficient services throughout the early election. However, it results in more than 25% of polling locations with waiting lines of more than 150 minutes. During the early election, our dynamic plan allows most polling locations to have waiting times less than 30 minutes, and around 1% of polling locations will have waiting times above 30 minutes. On election day, the dynamic plan allows approximately 50% of voting locations' voters to wait less than 30 minutes, more than 75% to wait less than 45 minutes, and 100% to wait less than 90 minutes. Due to fewer polling locations with long waiting times, the dynamic allocation is fairer than the fixed allocation.

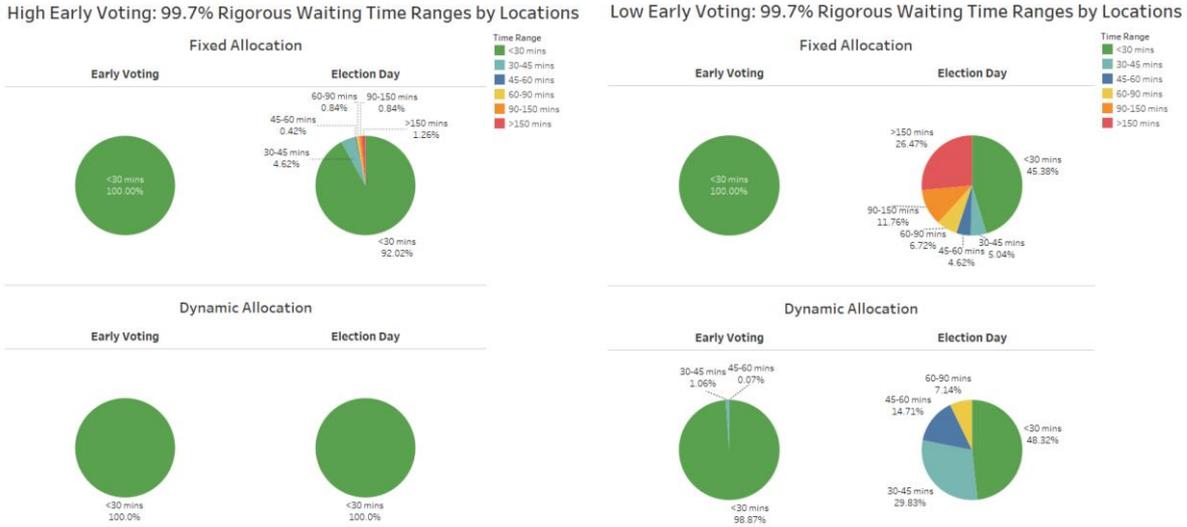

*Figure 5: Comparing Fixed Allocation and Dynamic Allocation's Wait Times*

In both scenarios, dynamic resource allocation requires fewer resources for early election and election day (Figure 6).

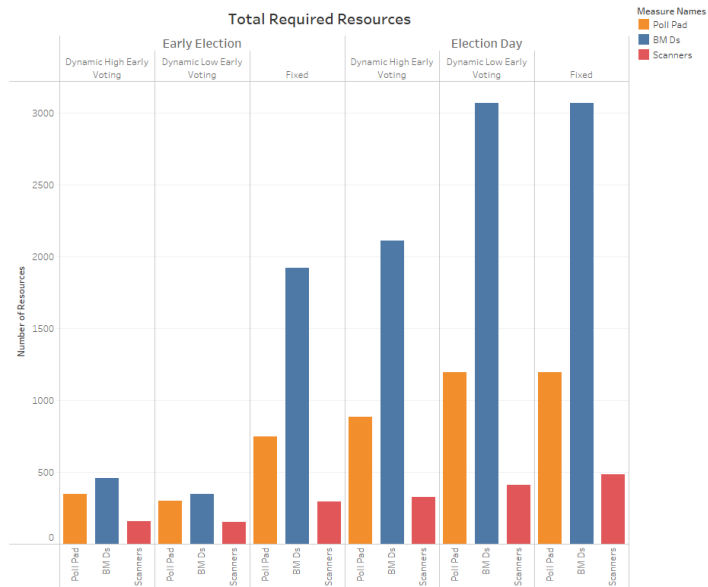

*Figure 6: Comparing Fixed Allocation and Dynamic Allocation's Required Resources*

Furthermore, dynamic resource allocation outperforms its fixed counterparts on average utilization rates of all three allocated resources in both scenarios (Table 1 and Table 2). This difference indicates that our dynamic methodology provides plans that exploit each piece of resource more scientifically, critically, and efficiently.

*Table 1: Fixed Allocation and Dynamic Allocation's Average Utilization in High Early Voting*

|  | *Early Voting* | | | *Election Day* | | |
|---|---|---|---|---|---|---|
|  | *Poll Pads* | *BMDs* | *Scanners* | *Poll Pads* | *BMDs* | *Scanners* |
| *Fixed* | 10.2% | 11.1% | 6.5% | 44.4% | 47.6% | 28.6% |
| *Dynamic* | 24.8% | 51.1% | 12.7% | 60.3% | 68.2% | 43.0% |

Table 2: Fixed Allocation and Dynamic Allocation's Average Utilization in Low Early Voting

|  | Early Voting | | | Election Day | | |
|---|---|---|---|---|---|---|
|  | *Poll Pads* | *BMDs* | *Scanners* | *Poll Pads* | *BMDs* | *Scanners* |
| *Fixed* | 6.4% | 6.9% | 4.1% | 68.2% | 73.0% | 43.9% |
| *Dynamic* | 16.3% | 39.5% | 8.0% | 72.0% | 75.2% | 53.2% |

We also want to point out that the fixed allocation considers social distancing, which puts a stronger location layout constraint on all polling locations (Barenji et al., 2023). Our model's allocation plan might require some polling locations to expand in practice. Policymakers should identify and improve such locations for better outcomes. Table 3 exemplifies a selected list of polling locations where our dynamic allocation improves performances and requires more resources than the fixed counterpart on election day in the low early voting scenario.

Table 3: Examples of Locations that Require More Resources on Election Day in Low Early Voting

|  | Fixed | | | | Dynamic | | | |
|---|---|---|---|---|---|---|---|---|
|  | *Poll Pads* | *BMDs* | *Scanners* | *Est. WT (mins)* | *Poll Pads* | *BMDs* | *Scanners* | *Est. WT (mins)* |
| *Abernathy Arts Center* | 5 | 11 | 2 | 90-150 | 6 | 17 | 2 | <30 |
| *Collier Park RC* | 3 | 6 | 1 | 45-60 | 3 | 9 | 1 | <30 |
| *Johns Creek High School* | 5 | 14 | 2 | 60-90 | 6 | 16 | 3 | <30 |
| *Morningside ES* | 6 | 18 | 3 | 60-90 | 9 | 19 | 4 | <30 |

## 6   Conclusion

This paper leverages PI principles and proposes a SEPI and a smart dynamic resource allocation methodology to ensure efficiency, fairness, resilience, and security in election systems. Through a case study of two scenarios, the dynamic methodology outperforms the fixed allocation.

This study has four main contributions. First, it incorporates queueing network models and lexicographic optimizations to test and optimize resource combinations of each polling location on each day. Second, the proposed resource allocation methodology requires fewer input resources and performs better than a traditional plan. Countries or regions with constrained resources can exploit it to allocate resources efficiently and critically. Third, this study identifies critical polling locations requiring more resources or expansion for polling policymakers. Fourth, this study promotes a scientific procedure in political decision-making. By adopting and applying PI principles and ideologies, we can ensure a safer and more effective outcome that protects equality and democracy around the globe.

This study also opens three future research avenues. The first is to effectively determine a limited number of polling locations with resource allocation plans, especially in regions with less compact populations. The second is to predict voter turnouts considering multiple factors, such as weather and media. The third is to allocate polling location workers during days.

Efficient, Fast, and Fair Voting Through Dynamic Resource Allocation
in a Secure Election Physical IntranetThis study finds dynamic resource allocations beneficial and critical for election systems in providing voters with satisfactory services. Policymakers and election practitioners are encouraged to utilize the proposed methodology.